\documentclass[12pt,reqno]{amsart}
\usepackage{amsfonts,amsmath,amscd,amssymb,amsthm}
\usepackage{caption,enumerate,graphicx,hyperref,perpage,url}
\usepackage[margin=2.6cm]{geometry}

\theoremstyle{plain}
\newtheorem{thm}{Theorem}

\newtheorem{cor}[thm]{Corollary}

\theoremstyle{remark}
\newtheorem*{defn}{\textbf{Definition}}
\newtheorem*{ex}{\textbf{Example}}
\newtheorem*{rmk}{\textbf{Remark}}

\numberwithin{equation}{section}

\MakePerPage{footnote} \allowdisplaybreaks 

\newcommand{\M}{\mathbb M}

\newcommand{\N}{\mathbb N}

\newcommand\R{\mathbb R}
\newcommand{\s}{\mathbb S}

\newcommand{\SL}{\mathrm{SL}}

\newcommand{\T}{\mathbb T}

\newcommand{\ve}{\varepsilon}

\newcommand{\Vol}{\mathrm{Vol}}

\newcommand{\Z}{\mathbb{Z}}

\title{From nodal points to non-equidistribution at the Planck scale}

\author{Xiaolong Han}
\email{xiaolong.han@csun.edu}
\address{Department of Mathematics, California State University, Northridge, CA 91330, USA}

\thanks{} 
\date{}

\begin{document}
\maketitle

\begin{abstract}
In this note, we make an observation that Laplacian eigenfunctions fail equidistribution at the Planck scale. Furthermore, equidistribution at the same scale also fails around the points where the eigenfunctions have large values.

Dans cette note, on fait une observation que les fonctions propres du laplacien ne s'equidistriburent \`a l'\'echelle de Planck. De plus, l'\'equidistribution \`a la m\^eme \'echelle ne sont plus valable au tour des points o\`u les functions propres ont des grand valeurs.
\end{abstract}

\section{Introduction}
Let $\M$ be a $n$-$\dim$ compact Riemannian manifold. Denote by $\Delta$ the non-negative Laplacian on $\M$. (If $\M$ has boundary, then we impose the Dirichlet or Neumann boundary condition.) We call $r:(0,\infty)\to(0,1)$ a \textit{small scale function} if $r(\lambda)\to0$ as $\lambda\to\infty$. We are concerned with the (non-)equidistribution of Laplacian eigenfunctions at small scales. 
\begin{defn}[Equidistribution at small scales]
Suppose that $r=r(\lambda)$ is a small scale function. Let $\{u_k\}_{k=1}^\infty$ be a sequence of eigenfunctions, $\Delta u_k=\lambda_k^2u_k$ with $\lambda_k\to\infty$ and $\|u_k\|_{L^2(\M)}=1$. We say that $\{u_k\}_{k=1}^\infty$ tend to be equidistributed at the scale $r$, if
\begin{equation}\label{eq:sse}
\int_{B(p,r_k)}|u_k|^2\,d\Vol=\frac{\Vol(B(p,r_k))}{\Vol(\M)}+o(r_k^n)\quad\text{as }k\to\infty,
\end{equation}
for all $p\in\M$ uniformly. Here, $r_k=r(\lambda_k)$, $B(p,r)\subset\M$ is the geodesic ball with center $p$ and radius $r$, and $d\Vol$ denotes the Riemannian volume on $\M$.
\end{defn}

The phenomenon of equidistribution of eigenfunctions has been extensively studied in physics and mathematics, see, e.g., Gutzwiller \cite[Chapter 15]{G} and Zelditch \cite[Chapter 9]{Z2}. In particular, Shnirelman \cite{Sh}, Zelditch \cite{Z1}, and Colin de Verdi\`ere \cite{CdV} proved that if the geodesic flow on a manifold $\M$ is ergodic, then any orthonormal basis of eigenfunctions (i.e., eigenbasis) contains a full density subsequence which tend to be equidistributed at the macroscopic scale, as a consequence, \eqref{eq:sse} holds with $r$ independent of $\lambda$. In fact, equidistribution of the subsequence holds in the phase space (a stronger condition than the one in the physical space $\M$). This result is called \textit{Quantum Ergodicity}, since  eigenfunctions describe the stationary states in the corresponding quantum system of the geodesic flow. 

It is well known that on a manifold with negative sectional curvature, the geodesic flow is ergodic, c.f., Katok-Hasselblatt \cite[Chapter 17]{KH}, so Quantum Ergodicity holds. On these negatively curved manifolds, Han \cite{Ha} and Hezari-Rivi\`ere \cite{HR} proved a weaker form of equidistribution at the scale $r(\lambda)=(\log\lambda)^{-\alpha}$ with some $\alpha>0$, that is, any eigenbasis contains a full density subsequence of eigenfunctions for which the two sides of \eqref{eq:sse} are uniformly comparable. Such small scale results have applications to problems including the $L^p$ norm estimates of eigenfunctions \cite{HR, So}. We refer to Zelditch's survey \cite{Z3} for the recent development of equidistribution at small scales.

In this note, we study scales $r$ at which equidistribution fails for all (real-valued) eigenfunctions (i.e., non-equidistribution) in the sense that one of the following two cases happens.
\begin{itemize}
\item Case I. For some $p\in\M$,
$$\int_{B(p,r)}|u|^2\,d\Vol\ll\Vol(B(p,r)).$$
\item Case II. For some $p\in\M$,
$$\int_{B(p,r)}|u|^2\,d\Vol\gg\Vol(B(p,r)).$$
\end{itemize}
We will make the conditions ``$\ll,\gg$'' precise in the theorems below. We shall also remark that the above inequalities can not hold for all points $p\in\M$ uniformly, due to the normalization that $\|u\|_{L^2(\M)}=1$. Moreover, the points at which the inequalities hold depend on the function $u$. 

Notice that an eigenfunction with eigenvalue $\lambda^2$ oscillates at a typical wavelength $\approx\lambda^{-1}$, which is usually referred as the \textit{Planck scale}. For example, the following picture demonstrates the density distribution of the eigenfunctions in the Barnett domain. As the eigenvalues increase from left to right, the oscillating wavelengths (Planck scale) decrease. 
\begin{figure}[h!]
\includegraphics[width=\textwidth]{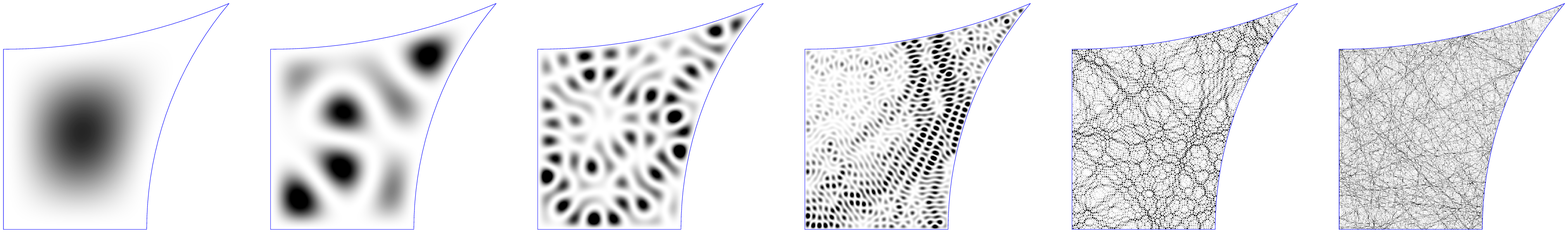}
\caption*{Credit: Alex Barnett}
\end{figure}

Concerning the (non-)equidistribution at small scales, we consider the  eigenfunctions on the unit circle $\R/2\pi\Z$ as the simplest model. In this case, the eigenvalues are $k^2$, $k\in\N$. The eigenfunctions are given by linear combinations of $\sin(kx)$ and $\cos(kx)$. These functions oscillate at the Planck scale $k^{-1}$ and clearly fail equidistribution at such scale. 

The goal of this note is to generalize the non-equidistribution result on the circle at the Planck scale to \textit{all} eigenfunctions on \textit{any} manifold. Indeed, we prove the non-equidistribution in Case I of eigenfunctions $u$ on balls at the Planck scale. The proof is based on two facts about the eigenfunctions:
\begin{enumerate}[A.]
\item The nodal points of $u$, $\Delta u=\lambda^2u$, are ``$\lambda^{-1}$ dense'' on $\M$. This fact is standard. See e.g., Colding-Minicozzi \cite[Lemma 1]{CM}: There is $a>0$ depending only on $\M$ such that for 
\begin{equation}\label{eq:r0}
R=a\lambda^{-1}
\end{equation} 
and any $p\in\M$, there is $q\in B(p,R/3)$ with $u(q)=0$. We fix $a$ throughout the note.
\item The average of the gradient of $u$, $|\nabla u|$, is bounded above by $O(\lambda)$ uniformly on $\M$. This fact follows from  Green's identity:
\begin{equation}\label{eq:av}
\int_\M|\nabla u|^2\,d\Vol=\int_\M u\Delta u\,d\Vol=\lambda^2\int_\M|u|^2\,d\Vol=\lambda^2,
\end{equation}
since $\|u\|_{L^2(\M)}=1$. By Cauchy-Schwarz inequality,
$$\int_\M|\nabla u|\,d\Vol\le\Vol(\M)^\frac12\left(\int_\M|\nabla u|^2\,d\Vol\right)^\frac12=\Vol(\M)^\frac12\lambda.$$ 
\end{enumerate}

We say that a collection of balls $\{B(p_j,r)\}_{j=1}^J$ is maximal disjoint in $\M$ if they are disjoint and any ball of radius $r$ in $\M$ intersects $B(p_j,r)$ for some $j=1,...,J$. Let $\{B(p_j,R)\}_{j=1}^J$ be a maximal disjoint collection of balls, in which $R=a\lambda^{-1}$. Then Fact A asserts that there is always a nodal point $q_j\in B(p_j,R/3)$, while Fact B implies that $|\nabla u|\lesssim\lambda$ if we average over all balls. Therefore, for ``almost all'' balls, $u$ cannot grow too fast from the nodal points to neighboring points. It then results non-equidistribution in Case I at the Planck scale on these balls. This is our main theorem:
\begin{thm}\label{thm:smallmass}
Let $R=a\lambda^{-1}$ and $\ve>0$. Then there exists $c>0$ depending on $\ve$ and $\M$ such that the following holds. Assume that $\Delta u=\lambda^2u$ and $\{B(p_j,R)\}_{j=1}^J$ is a maximal disjoint collection of balls in $\M$. Let $\delta=c\ve$ and $r\le\delta\lambda^{-1}$. Then there is a subcollection $\{p_{j_k}\}_{k=1}^K\subset\{p_j\}_{j=1}^J$ which satisfies 
\begin{enumerate}[(i).]
\item $K\ge(1-\ve)J$,
\item for all $k=1,...,K$, there is $q_{j_k}\in B(p_{j_k},R/3)$ such that
$$\int_{B(q_{j_k},r)}|u|^2\,d\Vol\le\ve\cdot\Vol(B(q_{j_k},r)).$$
\end{enumerate}
\end{thm}
Therefore, equidistribution of eigenfunctions fails on the balls $B(q_{j_k},r)$, in which $r\le\delta\lambda^{-1}$. These balls can be regarded as small Planck balls, each of which is contained in a large Planck ball $B(p_{j_k},R)$ with $R=a\lambda^{-1}$.

From Theorem \ref{thm:smallmass}, we have that the following corollary.
\begin{cor}\label{cor:smallmass}
Let $R=a\lambda^{-1}$ and $r=o(\lambda^{-1})$ be a small scale function. Assume that $\Delta u=\lambda^2u$ and $\{B(p_j,R)\}_{j=1}^J$ is a maximal disjoint collection of balls in $\M$. Then there is a subcollection $\{p_{j_k}\}_{k=1}^K\subset\{p_j\}_{j=1}^J$ which satisfies 
\begin{enumerate}[(i).]
\item $$\lim_{\lambda\to\infty}\frac KJ=1,$$
\item there is $q_{j_k}\in B(p_{j_k},R/3)$ such that
$$\lim_{\lambda\to\infty}\frac{\int_{B(q_{j_k},r)}|u|^2\,d\Vol}{\Vol(B(q_{j_k},r))}=0\quad\text{uniformly for all }k=1,...,K.$$
\end{enumerate}
\end{cor}

\begin{rmk}\hfill
\begin{itemize}
\item In Theorem \ref{thm:smallmass} and Corollary \ref{cor:smallmass}, the quantities $R$, $r$, $J$, and $K$ depend on the eigenvalue $\lambda^2$, whereas the collection of balls $\{B(q_{j_k},r)\}_{k=1}^K$ (on which equidistribution of the eigenfunction $u$ fails) depend on $u$.

\item According to (i) in Corollary \ref{cor:smallmass}, ``almost all'' the large Planck balls contains small Planck balls on which equidistribution fails. We shall mention that ``almost all'' can not be replaced by ``all''. See the example on the sphere in Section \ref{sec:nodal}. 
\end{itemize}
\end{rmk}

We next establish the non-equidistribution in Case II at the Planck scale around the points where the eigenfunction takes large values.

\begin{thm}\label{thm:smallvol}
There exists $\gamma>0$ depending on $\M$ such that the following holds. Given any $M>0$ and $p\in\M$ for which $|u(p)|\ge M$, we have that
$$\frac{\int_{B(p,r)}|u|^2\,d\Vol}{\Vol(B(p,r))}\ge\gamma M^2\quad\text{with }r=\gamma\lambda^{-1}.$$
Here, $\Delta u=\lambda^2u$ and $M$ can depend on $\lambda$.
\end{thm}
\begin{rmk}\hfill
\begin{itemize}
\item Suppose that there is a sequence of eigenfunctions $\{u_k\}_{k=1}^\infty$ for which $\|u_k\|_{L^\infty(\M)}\to\infty$ as $k\to\infty$. Then as an immediate consequence of Theorem \ref{thm:smallvol}, non-equidistribution at the Planck scale in Case II happens around the points where $u_k$ achieves the maximum value.
\item The eigenfunctions on rectangles with irrational ratio are all sine or cosine functions (or finite linear combinations of these functions) so are uniformly bounded. It is plausible that they are the only manifolds on which all eigenfunctions are uniformly bounded. See Toth-Zelditch \cite{TZ} for related discussion. But this is unknown to the author's knowledge. 
\end{itemize}
\end{rmk}

Moreover, we prove non-equidistribution at larger scales than the Planck one, depending on $M$ in Theorem \ref{thm:smallvol}.
\begin{thm}\label{thm:smallvolgen}
Let $\ve>0$. Then there exists $\gamma>0$ depending on $\ve$ and $\M$ such that the following holds. Given any $M\ge1$ and $p\in\M$ for which $|u(p)|\ge M$, we have that
$$\frac{\int_{B(p,r)}|u|^2\,d\Vol}{\Vol(B(p,r))}\ge\frac1\ve\quad\text{with }r=\gamma\lambda^{-1}M^{\frac 2n}.$$
Here, $\Delta u=\lambda^2u$ and $M$ can depend on $\lambda$.
\end{thm}
\begin{rmk}\hfill
\begin{itemize}
\item From the Weyl-type estimate that $\|u\|_{L^\infty(\M)}\le C\lambda^{(n-1)/2}$ by H\"ormander \cite{Ho}, we see that the best scale that one can hope to get from Theorem \ref{thm:smallvolgen} is $\lambda^{-1}(\lambda^{(n-1)/2})^{2/n}=\lambda^{-1/n}$ (which is much larger than the Planck scale $\lambda^{-1}$.)
\item On the spheres $\s^n$, there are Laplacian eigenfunctions that saturate the above $L^\infty$ bound, e.g., the zonal harmonics. Therefore, these eigenfunctions display non-equidistribution in Case II at the scale $\lambda^{-1/n}$. It is unclear whether the sphere is the only example of manifolds for non-equidistribution at such scale to hold.
\end{itemize}
\end{rmk}

\subsection*{Related literature}
The results in this note are on general manifolds. If there is additional arithmetic structure, then some non-equidistribution results at various scales are known. 

On the tori $\T^n=\R^n/2\pi\Z^n$, Bourgain proved non-equidistribution in Case I for eigenfunctions at scales $r$ such that $r=o(\lambda^{-1/(n-1)})$, which is much larger than the Planck one. (The result was published in \cite[Theorem 4.1]{LR}.) However, the non-equidistribution in \cite{LR} is for balls with a fixed center, as opposite to balls that are separated by the Planck scale in Theorem \ref{thm:smallmass} and Corollary \ref{cor:smallmass}. In addition, the method used in \cite{LR} differ with the one in this note. That is, the eigenfunction $u$, $\Delta u=\lambda^2u$, on $\T^n$ can be written as
$$u(x)=c\sum_{k\in\Z^n\text{ with }|k|=\lambda}a_k\sin(k\cdot x)+b_k\cos(k\cdot x),$$  
in which $a_k,b_k\in\R$ and $c\in\R$ is the normalizing factor such that $\|u\|_{L^2(\T^d)}=1$. Bourgain selected a special collection of lattice points $k$ on the sphere $S_\lambda=\{x\in\R^n:|x|=\lambda\}$. Cancellation can then be exploited so that $\int_{B(p,r)}|u|^2/\Vol(B(p,r))\to0$ for $r=o(\lambda^{-1/(n-1)})$ and some $p\in\T^n$. On a general manifold, the arithmatic structure is not available. Our proof of Theorem \ref{thm:smallmass} and Corollary \ref{cor:smallmass} rely on the elliptic estimates of eigenfunctions at the Planck scale.

On the modular surface $\SL(2,\R)/\SL(2,\Z)$, Humphries \cite[Theorem 1.14]{Hu} proved non-equidistribution in Case II for Hecke-Maass eigenforms (Laplacian eigenfunctions with additional arithmetic structure) at scales $\lambda^{-1}(\log\lambda)^\alpha$ ($\alpha>0$) around the Heegner points. The proof is similar to the one of Theorem \ref{thm:smallvolgen}. That is, the eigenfunction $u$ achieves values of order $(\log\lambda)^\alpha$ at the Heegner points, see Mili\'cevi\'c \cite{M}. Non-equidistribution at the scale $\lambda^{-1}(\log\lambda)^\alpha$ therefore follows from Theorem \ref{thm:smallvolgen} by setting $M=(\log\lambda)^\alpha$.

\section{Non-equidistribution around nodal points}\label{sec:nodal}
We prove Theorem \ref{thm:smallmass} in this section. Let $\Delta u=\lambda^2u$ and $R=a\lambda^{-1}$ in \eqref{eq:r0}. Assume that $\{B(p_j,R)\}_{j=1}^J$ is a maximal disjoint collection of balls in $\M$. Then $J\ge c_0\lambda^n$ for some $c_0>0$ that depends only on $\M$.

According to the mean value inequality for $\nabla u$ (see Schoen-Yau \cite[Section II.6]{SY} and Zelditch \cite[Section 5.3.4]{Z2}), we have that
$$\sup_{B\left(p_j,\frac{2}{3}R\right)}\{|\nabla u|^2\}\le\frac{C_0}{\Vol(B(p_j,R))}\int_{B\left(p_j,R\right)}|\nabla u|^2\,d\Vol\le C_1\lambda^n\int_{B\left(p_j,R\right)}|\nabla u|^2\,d\Vol,$$
in which $C_1>0$ depends only on $\M$. It then follows that
\begin{eqnarray*}
\frac1J\sum_{j=1}^J\sup_{B\left(p_j,\frac{2}{3}R\right)}\{|\nabla u|^2\}&\le&\frac{1}{c_0\lambda^n}\sum_{j=1}^JC_1\lambda^n\int_{B\left(p_j,R\right)}|\nabla u|^2\,d\Vol\\
&\le&\frac{C_1}{c_0}\sum_{j=1}^J\int_{B(p_j,R)}|\nabla u|^2\,d\Vol\\
&\le&\frac{C_1}{c_0}\int_\M|\nabla u|^2\,d\Vol\\
&=&C_2\lambda^2,
\end{eqnarray*}
in which $C_2=C_1/c_0$ depends only on $\M$ and the last step follows from \eqref{eq:av}. 

By the Tchebychev inequality,
$$\frac1J\cdot\#\left\{j=1,...,J:\sup_{B\left(p_j,\frac{2}{3}R\right)}\{|\nabla u|^2\}>\frac{C_2\lambda^2}{\ve}\right\}\le\ve$$
Denote $\{p_{j_k}\}_{k=1}^K\subset\{p_j\}_{j=1}^J$ the subcollection that
$$\sup_{B\left(p_{j_k},\frac23R\right)}|\nabla u|^2\le\frac{C_2\lambda^2}{\ve}.$$
Then we have that $K\ge(1-\ve)J$, proving Condition (i) in Theorem \ref{thm:smallmass}.

Now for each $k=1,...,K$, there is $q_{j_k}\in B(p_{j_k},R/3)$ such that $u(q_{j_k})=0$ by \eqref{eq:r0}. For $\delta<a/3$  and $r\le\delta\lambda^{-1}\le R/3$, $B(q_{j_k},r)\subset B(p_{j_k},2R/3)$. Let $d(q_1,q_2)$ denote the Riemannian distance between $q_1,q_2\in\M$. 

For any $q\in B(q_{j_k},r)$, 
\begin{equation}\label{eq:mvt}
|u(q)-u(q_{j_k})|\le\sup_{B(q_{j_k},r)}\{|\nabla u|\}\cdot d(q,q_{j_k}).
\end{equation}
Since $u(q_{j_k})=0$ and $B(q_{j_k},r)\subset B(p_{j_k},2R/3)$,
\begin{eqnarray*}
\int_{B(q_{j_k},r)}|u(q)|^2\,d\Vol&\le&\sup_{B(q_{j_k},r)}\{|\nabla u|^2\}\cdot\int_{B(q_{j_k},r)}d(q,q_{j_k})^2\,d\Vol\\
&\le&\sup_{B\left(p_{j_k},\frac{2}{3}R\right)}\{|\nabla u|^2\}\cdot\int_{B(q_{j_k},r)}d(q,q_{j_k})^2\,d\Vol\\
&\le&\frac{C_2\lambda^2}{\ve}\cdot r^2\Vol(B(q_{j_k},r)),
\end{eqnarray*}
Therefore, for any $\ve>0$,
$$\int_{B(q_{j_k},r)}|u|^2\,d\Vol\le\frac{C_2\lambda^2r^2}{\ve}\cdot\Vol(B(q_{j_k},r))\le\ve\Vol(B(q_{j_k},r))$$
for all $r\le\delta\lambda^{-1}$ with $\delta=\ve/\sqrt{C_2}$.

\begin{rmk}
The above proof uses the simple observation that the non-equidistribution in Case I around nodal points is a direct consequence of the gradient bound. We use the average gradient bound that $|\nabla u|\lesssim\lambda$ to prove non-equidistribution at the Planck scale $\lambda^{-1}$ at ``almost all'' balls in Theorem \ref{thm:smallmass} and Corollary \ref{cor:smallmass}. 

On the other hand, recall the Weyl-type estimate 
$$\left\|\nabla u\right\|_{L^\infty(\M)}\lesssim\lambda^{\frac{n+1}{2}}.$$
See e.g., Sogge-Zelditch \cite{SZ}. Therefore, we immediately have non-equidistribution of $u$ at the scale $\lambda^{-(n+1)/2}$ around \textit{every} nodal point. But notice that this scale is much smaller than the Planck one.
\end{rmk}

\begin{ex}[Highest weight spherical harmonics]
We shall point out that the non-equidistribution in Case I may not hold everywhere, that is, the ``almost all'' condition in Theorem \ref{thm:smallmass} and Corollary \ref{cor:smallmass} is optimal. For example, in the polar coordinates $\theta\in[0,\pi]$ and $\phi\in[0,2\pi)$ on the sphere $\s^2$, let $u_k(\theta,\phi)=k^{1/4}(\sin\theta)^k\sin(k\phi)$. Then $u_k$ is a Laplacian eigenfunction on $\s^2$ with eigenvalue $k(k+1)$. The factor $k^{1/4}$ normalizes the function so that $c_1\le\|u_k\|_{L^2(\s^2)}\le c_2$, in which the positive constants $c_1$ and $c_2$ are absolute. See, e.g., Zelditch \cite[Section 4.4.5]{Z2}. 

The function $u_k$ has mass concentrated in an $O(k^{-1/2})$ neighborhood of the equator $\{(x,y,z)\in\s^2:z=0\}$ and are commonly referred as the ``highest weight spherical harmonics''. In particular, 
$$(\sin\theta)^k=\left|\cos\left(\frac\pi2-\theta\right)\right|^k\ge\frac12\quad\text{if }\left|\theta-\frac\pi2\right|\le\frac{1}{10}k^{-\frac12}.$$
Fix any $\delta\le1$ and let $r=\delta k^{-1}$. Take any ball $B(q,r)\subset\s^2$ such that $q=(\theta_q,\phi_q)$ with $|\theta_q-\pi/2|<k^{-1/2}/20$. Then it is straightforward to compute that
$$\int_{B(q,r)}|u_k|^2\,d\Vol\ge ck^\frac12\cdot\Vol(B(q,r)),$$
in which the positive constant $c$ depends only on $\delta$. This shows that the non-equidistribution in Case I at the Planck scale can not happen for balls in the strip. However, for any maximal disjoint balls with separation $\approx k^{-1}$, the ones in this strip are of ``zero density'' as $k\to\infty$, which agrees with Condition (i) in Corollary \ref{cor:smallmass}.
\end{ex}

\begin{rmk}
We use the same example as above to show that the scale of non-equidistribution in Case I at \textit{single} points can be much larger than the Planck one. 

Take $p=(\theta_p,\phi_p)$ with $\theta_p=0$ or $\pi$ (i.e., $p$ is the north or south pole) and $r<\pi/2$. So $u_k(p)=0$ and $B(p,r)$ does not intersect the equator. Then
$$\int_{B(p,r)}|u_k|^2\,d\Vol\le Ck^{\frac12}\int_0^{\cos r}t^kt\,dt\le Ck^{-\frac12}(\cos r)^{k+2}.$$
Hence,
$$\frac{\int_{B(p,r)}|u_k|^2\,d\Vol}{\Vol(B(p,r))}\le\frac{Ck^{-\frac12}(\cos r)^{k+2}}{r^2}\to0\quad\text{as }k\to\infty.$$
That is, one observes non-equidistribution of $u_k$ in Case I at any scale $r<\pi/2$ around the north and south poles. This is a reflection of the fact that $u_k$ and $\nabla u_k$ decay exponentially fast away from the equator.
\end{rmk}

\section{Non-equidistribution around points with large values}
We prove Theorem \ref{thm:smallvol} in this section. For $\delta>0$ and $r=\delta\lambda^{-1}$, the mean value inequality yields that
$$\sup_{B\left(p,\frac r2\right)}\{|u|^2\}\le\frac{c}{\Vol(B(p,r))}\int_{B(p,r)}|u|^2\,d\Vol.$$
Here, $c>0$ depends only on $\M$. See Schoen-Yau \cite[Section II.6]{SY} and Zelditch \cite[Corollary 5.11]{Z2} for the mean value inequality applied to eigenfunction in such balls. Thus,
$$\int_{B(p,r)}|u|^2\,d\Vol\ge c\Vol(B(p,r))\cdot\sup_{B\left(p,\frac r2\right)}\{|u|^2\}\ge c\Vol(B(p,r))M^2,$$
because $|u(p)|\ge M$. Therefore, 
$$\frac{\int_{B(p,r)}|u(y)|^2\,d\Vol}{\Vol(B(p,r))}\ge cM^2\ge\delta M^2,$$
by choosing $\delta\le c$.

Theorem \ref{thm:smallvolgen} follows by noticing that if $r=\delta\lambda^{-1}M^{2/n}$, then
$$\int_{B(p,r)}|u|^2\,d\Vol\ge\int_{B\left(p,\delta\lambda^{-1}\right)}|u|^2\,d\Vol.$$

\section*{Acknowledgments}
I learned the problem to prove non-equidistribution of Laplacian eigenfunctions at the Planck scale from Peter Humphries in June 2020. I want to thank him and Melissa Tacy for the discussions about this problem and Stephane Nonnenmacher for the suggestion to use the ``average'' gradient bound in Theorem \ref{thm:smallmass}. I also want to thank Alex Barnett for very helpful comments and the permission to use the picture he generated in this note.

\end{document}